\documentclass[a4,11pt]{article}
\usepackage[sumlimits]{amsmath}
\usepackage{amssymb,amscd}
 \usepackage[francais]{babel}

\usepackage{graphicx}  

\def\nal{\| \hspace{-.10em} |}

\def\1{{\mathbf 1 }}

\catcode`\é=\active\def é{\'e} \catcode`\è=\active\def è{\`e}
\catcode`\ç=\active\def ç{\c{c}} \catcode`\à=\active\def à{\`a}
\catcode`\ê=\active\def ê{\^e} \catcode`\ù=\active\def ù{\`u}
\catcode`\ô=\active\def ô{\^o} \catcode`\û=\active\def û{\^u}
\catcode`\î=\active\def î{\^{\i}} \catcode`\â=\active\def â{\^a}
\catcode`\ö=\active\def ö{\"o}
\catcode`\ë=\active\def ë{\"e}

 \def\ooo{{\cal O}}
  
\def\sss{{\cal S}}

\def\ph{\varphi}
\def\eps{{\varepsilon}}

\def\Om{{\Omega}}

\def\VaO{V_{\alpha}(\Omega)}

\def \Nal{\mathcal N}

\def \Ual{\mathcal U}
\def \Val{\mathcal V}

\newtheorem{theo}{Theorem}
\newtheorem{lem}[theo]{Lemma}

\newtheorem{prop}[theo]{Proposition}
\def\N{{\mathbb N}}
\def\R{{\mathbb R}}
\def\Z{{\mathbb Z}}

\def\Om{{\Omega}}

\def\ph{\varphi}
\def\eps{{\varepsilon}}

\def\cov{{\rm Cov}}

\def\1{{\mathbf 1}}

\title{Exponential decay of correlations for a real-valued
dynamical system generated by a $k$ dimensional system}

\author{JAGER Lisette, MAES Jules, NINET Alain}

\begin{document}

\maketitle

\section*{Abstract}

\noindent We study the real, bounded-variables process $\{X_n, n\in \N \}$ 
 defined by a $k$-term recurrence relation
 $X_{n+k} = \varphi(X_n, \hdots, X_{n+k-1})$. We prove the
 decay of correlations, mainly under purely analytic hypotheses concerning the 
function $\varphi$ and its partial derivatives.

%
%

\section{Introduction}

\noindent
Since the eighties, many statisticians have studied nonlinear time series in order
 to model various phenomena in Physics, Economics and Finance. About this
 subject one may consult, for example,  Chan-Tong \cite{TON1},
Tong \cite{TON2} and Gu\'egan \cite{GU}. With chaos theory it became clear that
 such a series is a perturbed dynamical system. 
 For an extensive survey of litterature about chaos theory, one may read
  Collet-Eckmann \cite{CE}, Lasota-Mackey \cite{LM} or Liverani \cite{LIV}. To undertake the study of
nonlinear time series with the help of the theory of dynamical systems, 
as a first step is the study of a non-perturbed dynamical system
defined by a recurrence relation of order $k$  (which Tong names ``skeleton''
in \cite{TON2}).
Indeed, we consider the bounded variables model, 
$X_{n+k}= \ph(X_n, \hdots, X_{n+k-1})$, where $\varphi$ is piecewisely defined 
on the set   $\Ual^k$ and takes its values in $\Ual$,
 where $\Ual=[-L,L]$ for $L \in \R_+^*$. This model gives rise to
a dynamical system  $ (\Om,\tau ,\mu,T)$ where $\Om$ is a compact subset of 
 $\R^k$ and  $\mu$ is a measure preserved by the transformation
 $T :\Om \rightarrow \Om$.
Under conditions  on $\varphi$, which ensure that $T$ satisfies Saussol's  
hypotheses  \cite{SAU}, we obtain the exponential decay of correlations
if $T$ is mixing. More precisely, for well-chosen applications $f$
 and $g$, we prove that there exist  constants  $C=C(f,h)>0$ and $\rho\in ]0,1[ $ such that:
 $$ \left|  \int_{\Om} f \circ T^n \,  h \ d\mu - \int_{\Om} f d\mu ~ \int_{\Om}h d\mu  \right| \leqslant C \, \rho^n.  $$
This result can be stated in the following way, in the case when $X_0$ 
has distribution $\mu$:
$$ \left| \, \cov ( \, f(X_n),h(X_0) \,) \,  \right| \leqslant C \, \rho^n .$$ 
Other methods could give the same result, under different hypotheses on
the induced system. One can study recurrence times, like Young 
\cite{YOU} or Gou\"ezel \cite{GO}.
 For a general view of all these different techniques, see 
 Alves-Freitas-Luzzato-Vaienti \cite{AFLV}.\\
We had studied the case when $k=2$ in a precedent paper \cite{JMN}. It was then
 possible to give a very precise result about the localization of the
 eigenvalues of the $2\times 2$ matrix (\ref{CoeffB}), which expresses the
 fact that an auxiliary transformation is expanding. This allows us
to obtain better estimates than in the general case, which is the subject
 of the present work.\\
At the end of this article, we illustrate our results by applying them to
a nonlinear example.

%
%

\section{Hypotheses and results}

\noindent Let $L\in \R_+^*$ and let us consider an application
 $\varphi : [-L,L]^k \rightarrow [-L,L]$,
 defined piecewisely  on $[-L,L]^k$.\\
Under conjugation by an affine function, similar results could be obtained
for an application $\varphi$ defined on  $[a,b]^k$, with values in  $[a,b]$.\\

\noindent Suppose that all following conditions are fulfilled:
\begin{enumerate}
\item \label{hyp1}
There exists $d \in \N^*$ such that
$$
[-L,L]^k= \bigcup_{j=1}^d O_j \ \cup \Nal,
$$
where the  $O_j$ are nonempty open subsets, $\Nal$ is Lebesgue negligible and 
the union is disjoint. The boundary of each $O_j$ is contained in a compact, 
$C^1$, $(k-1)-$ dimensional submanifold of $\R^k$.
\item \label{hyp2}
There exists $\eps_1>0$ such that, for every $j\in \{1,\dots, d\}$, there
 exists a map  $\ph_j$ defined on   $B_{\eps_1}(\overline{O_j}) = \{ (x_1, \hdots, x_k) \in \R^k, ~ d((x_1, \hdots, ,x_k),\overline{O_j}) \leq \eps_1\}$ with values in  $\R$ and satisfying $\ph_j|_{O_j} = \ph|_{O_j}$.
\item \label{hyp3}
The application  $\ph_j$ is bounded and  $C^{1,\alpha}$ on 
$B_{\eps_1}(\overline{O_j})$  for an  $\alpha\in ]0,1]$ \footnote{If $\ph_j$ is   $C^2$ on $B_{\eps_1}(\overline{O_j})$, it is necessarily  $C^{1,\alpha}$ on $B_{\eps_1}(\overline{O_j})$ with $\alpha = 1$}, which means that  $\ph_j$ is  $C^1$
and that there exists  $C_j>0$ such that, for all $ (u_1, \hdots, u_k), (v_1, \hdots, v_k)$ in $B_{\eps_1}(\overline{O_j})$,  all  $i \in \{1, \dots, k\}$ :
$$
\left|\frac{\partial \ph_j}{\partial x_i}(u_1, \hdots, u_k) -
\frac{\partial \ph_j}{\partial x_i}(v_1, \hdots, v_k)\right|
\leq C_j ||(u_1, \hdots, u_k)-(v_1, \hdots, v_k)||^{\alpha}.
$$
We assume that there exist constants 
 $A>1$ and $\sigma >1$ satisfying $A^{2/k} > \sigma$, 
 $M\in ~ ]0,M_0(\sigma,A)[$ such that for all $2 \leq i \leq k$ :
$$ 
 \forall (u_1, \hdots, u_k) \in B_{\eps_1}(\overline{O_j}),\quad
 \left|\frac{\partial \ph_j}{\partial x_1}(u_1, \hdots, u_k)
\right|\geq A,\quad
\left|\frac{\partial \ph_j}{\partial x_1}(u_1, \hdots, u_k) \times \frac{\partial \ph_j}{\partial x_i}(u_1, \hdots, u_k) \right| \leq M
$$
with :
$$ M_0(\sigma,A) = \frac{-(k-1) \gamma^{k-1} + \sqrt{(k-1)^2 \gamma^{2k-2}+4(k-2)\gamma^{2k+1}(\frac{1}{\gamma^2}-\sigma)}}{2(k-2) \gamma^{2k+1}} >0 
$$
and
$$  \gamma = A^{-1/k}.$$
These very tight conditions are due to
the loss of precision in the localization of the eigenvalues of the matrix
 $B$ - see  (\ref{CoeffB}) - in the case when $k>2$.
\item \label{hyp4}
The sets $O_j$ satisfy the following geometrical condition: 
\footnote{In favorable cases, the geometrical hypothesis can be replaced by the following one, stronger but much simpler: for all points  $(u_1, u_2, \hdots, u_k)$ and $(v_1, u_2, \hdots, u_k)$ in $B_{\eps_1}(\overline{O_j})$, the segment $[(u_1, u_2, \hdots, u_k),(v_1, u_2, \hdots, u_k)]$ is contained in  $B_{\eps_1}(\overline{O_j})$} for all $(u_1, u_2, \hdots, u_k)$ and $(v_1, u_2, \hdots, u_k)$ in $ B_{\eps_1}(\overline{O_j})$, there exists a $C^1$ path  $\Gamma=(\Gamma_1,\hdots,\Gamma_k) : [0,1]\rightarrow  B_{\eps_1}(\overline{O_j})$ between  $(u_1, u_2, \hdots, u_k)$ and $(v_1, u_2, \hdots, u_k)$, with nonzero gradient  satisfying
\begin{equation*}
\forall ~ t\in ]0,1[, ~ \left|\Gamma_1'(t)\right|
 >\frac{M}{A^2} \sum_{i=2}^k \left|\Gamma_i'(t)\right| .
\end{equation*}
\item \label{hyp5}
The maximal number of $C^1$ arcs of  $\Nal$ crossing  is $Y \in \N^*$.
Moreover, one sets
$$
s = \frac{1}{\sqrt{\sigma}} < 1
$$
and imposes that
$$
\eta:=  s^{\alpha} + \frac{4 s}{1-s}Y \frac{\gamma_{k-1}}{\gamma_k} <1
$$
where  $\displaystyle \gamma_k = \frac{\pi^{k/2}}{ \Gamma(\frac{k}{2}+1)}$
is the volume of the unit sphere of $\R^k$.
\end{enumerate}
(This last condition is the most  restricting of all, for it gives an
 upper bound for $s$ and thus a lower bound for $\sigma$ and $A$.)\\
\noindent For every  $j \in \{1, ..., d\}$, one denotes  by $U_j$ (resp. 
$W_j$, $\Nal'$)  the image of $O_j$ (resp. $B_{\eps_1}(\overline{O_j})$, $\Nal$)
 under the transformation which associates 
with $(u_1, \hdots, u_k) \in \R^k$ the
 point $(u_1,\gamma u_2, \hdots, \gamma^{k-1} u_k)$. The set
 $\Om = [-L,L] \times [-\gamma L, \gamma L] \times \hdots \times [-\gamma^{k-1}L, \gamma^{k-1}L]$, with which we shall work, is the image of $[-L,L]^k$ under the same transformation.\\
For every non-negligible Borel set $S$ of
 $\R^k$, for every $f \in L^1_m(\R^k,\R)$, one sets :
$$ Osc(f,S) = \underset{S}{Esup} f - \underset{S}{Einf} f $$
where $\underset{S}{Esup}$ et $\underset{S}{Einf}$ are the essential supremum
 and infimum on $S$ with respect to the Lebesgue measure $m$.\\
One then defines the norm  $\| \cdot \|_{\alpha}$ by
$$
|f|_{\alpha}= \sup_{0<\eps<\eps_1}\eps^{-\alpha} \int_{\R^k}  {\rm Osc}(f,B_{\eps}(x_1, \hdots, x_k))\ dx_1 \hdots dx_k \ , \qquad \| f \|_{\alpha} = \| f \|_{L^1_m} + |f|_{\alpha}
$$
and the set
$V_{\alpha} = \{ f \in L^1_m(\R^k,\R), ~ \| f \|_{\alpha} < +\infty \}$.\\

\noindent
We introduce similar notions on  $\Om$ : for every $0 < \eps_0 < \gamma^{k-1} L$, for every $g \in L^{\infty}_m(\Om,\R)$, we define:
$$
N(g,\alpha,L) = \sup_{0<\eps<\eps_0}\eps^{-\alpha} \int_{\Om} {\rm Osc}(g,B_{\eps}(x_1, \hdots, x_k) \cap \Om)\ dx_1 \hdots dx_k.
$$
One then sets:
$$
||g||_{\alpha,L}=  N(g,\alpha,L) + 2 K(\Omega) ~ \eps_0^{1-\alpha} ||g||_{\infty} + ||g||_{L^1_m}
$$
where $K(\Omega) = 2^{k+2} (\sum\limits_{i=1}^k 2 \gamma^{i-1} L)^{k-1} = 2^{2k+1} L^{k-1} (\frac{1-\gamma^k}{1-\gamma})^{k-1}$.\\
The function  $g$ is said to belong to  $\VaO$ if this expression is finite. 
This set does not depend on the choice of $\eps_0$, but $N$ and $\| . \|_{\alpha,L}$ do.\\
There exist relations between the sets $V_{\alpha}(\Omega)$ and $V_{\alpha}$.
Indeed, one can prove the following result using Proposition 3.4 of \cite{SAU}:

\begin{prop}{\ }\\
\begin{enumerate}
\item
If $g \in \VaO$ and if one extends $g$ to a function $f$
defined on $\R^k$, setting  $f(x)=0$ if $x\notin \Om$, 
then $f \in V_{\alpha}$ and
$$
\| f \|_{\alpha} \leq \| g \|_{\alpha,L}.
$$
\item
Conversely
let  $f\in V_{\alpha}$ and set $g= f \1_{\Om}$. Then $g\in \VaO$ and the following holds:
$$
 \| g \|_{\alpha,L}\leq \left( 1+ 2 K(\Omega)
 \frac{\max(1,\eps_0^{\alpha})}{\gamma_k ~ \eps_0^{k-1+\alpha}}\right)
\| f \|_{\alpha}.
$$
\end{enumerate}
\end{prop}

\noindent Under the hypotheses 1-5 listed above, we  obtain a first result:

\begin{theo}\label{resconjug}
Let $T$ be the  transformation defined on  $\Om$ by : $\forall  u=(u_1, \dots, u_k) \in U_j$ :
$$
T(u) = T_j(u)= \left( \frac{u_2}{\gamma}, ~ \dots, \frac{u_k}{\gamma}, ~ \gamma^{k-1} \ph_j(u_1, \frac{u_2}{\gamma}, ~ \dots, ~ \frac{u_k}{\gamma^{k-1}})\right).
$$
The applications $T_j$ can be defined  naturally on $W_j$ by the same
 formula. Then 
\begin{enumerate}
\item 
The Frobenius-Perron operator $P : L^1_m(\Om) \rightarrow L^1_m(\Om)$ associated
 with $T$  has a finite number of eigenvalues of  modulus $1$, 
 $\lambda_1,\dots,\lambda_r$.
\item 
For every $i\in\{1,\dots,r\}$, the eigenspace  
$E_i=\{ f\in L^1_m(\Om) \ : \ Pf=\lambda_i f\}$  associated with the eigenvalue
 $\lambda_i$ is finite-dimensional and contained in $V_{\alpha}(\Om)$.
\item 
The operator $P$ decomposes as
$$
P= \sum_{i=1}^r \lambda_i P_i + Q,
$$
where the  $P_i$ are projections on the spaces $E_i$, $\nal P_i \nal_1\leq 1$
and  $Q$ is a linear operator defined on $ L^1_m(\Om)$, such that
  $Q(\VaO) \subset \VaO$, $\sup\limits_{n \in \N^*}\nal Q^n \nal_1<\infty$ and
  $\nal Q^n \nal_{\alpha,L} = O(q^n)$ when $n \rightarrow +\infty$, for a given
  $q\in ]0,1[$. Moreover, $P_iP_j=0$ if $i\neq j$, $P_iQ=QP_i=0$ for every $i$.
\item 
The operator $P$ has the eigenvalue $1$. Set $\lambda_1=1$, let 
 $h_*=P_1 \1_{\Om}$ and  $d\mu=h_* ~ dm$. Then $\mu$ is the greatest absolutely continuous invariant measure (ACIM) of $ T$, which means that, if
 $\nu <<m$ and if  $\nu$ is $T$-invariant, then $\nu<<\mu$.
\item
 The support of $\mu$ can be decomposed into a finite number of mutually 
disjoint measurable sets, on which a power of $T$ is mixing. More precisely,
 for every $j \in \{1,2,\dots, \dim(E_1)\}$, there exist a number
 $L_j \in \N^*$ and $L_j$ mutually disjoint sets $W_{j,l}$ $(0 \leq l \leq L_{j}-1)$, satisfying  $T(W_{j,l})=W_{j,l+1 \mod(L_j)}$, $T^{L_j}$ being mixing on every  $W_{j,l}$.
One denotes by  $\mu_{j,l}$ the normalized restriction of $\mu$ on  $W_{j,l}$,
defined by
$$
\mu_{j,l}(B)= \frac{\mu(B\cap W_{j,l})}{\mu(W_{j,l})}, \ d\mu_{j,l} = \frac{h^* \1_{W_{j,l}} }{\mu(W_{j,l})} dm.
$$
Saying that $T^{L_j}$ is mixing on every $W_{j,l}$ means that, for every $f \in L^1_{\mu_{j,l}}(W_{j,l})$ and every $h \in L^{\infty}_{\mu_{j,l}}(W_{j,l})$,
$$ \lim\limits_{n \rightarrow + \infty} <T^{nL_j} f,h>_{\mu_{j,l}} = <f,1>_{\mu_{j,l}} <1,h>_{\mu_{j,l}} $$
with indifferently used notations:
 $<f,g>_{\mu'} = \mu'(fg) = \int f g ~ d\mu'$. 
\item Moreover, there exist real constants $C>0$ and $0<\rho < 1$ such that, 
for every $h$ in $\VaO$ and $f\in L^1_{\mu}(\Om)$, the following holds:
$$
\left| \int_{\Om} f \circ T^{ n \times ppcm(L_i)}  h \ d\mu -\sum_{j=1}^{\dim(E_1)} \sum_{l=0}^{L_j-1} \mu(W_{j,l}) <f,1>_{\mu_{j,l}} <1,h>_{\mu_{j,l}} \right| \leq C  ||h||_{\alpha,\Om} ||f||_{L^1_{\mu}(\Om)}~\rho^{n }.
$$
\item If, moreover, $T$ is mixing,
 \footnote{Which is equivalent to: if $1$ is the only modulus-$1$ eigenvalue
 of $P$ and if, additionnaly, it is simple}
the preceding result can be  stated as: there exist real constants $C>0$ and
 $0<\rho < 1$  such that, for every $h$ in $\VaO$ and $f\in L^1_{\mu}(\Om)$,
 one has:
$$ \left|  \int_{\Om} f \circ T^n \,  h \ d\mu - \int_{\Om} f d\mu ~ \int_{\Om} h d\mu  \right| 
 \leq C  ||h||_{\alpha,\Om}~||f||_{L^1_{\mu}(\Om)}~\rho^{n }. $$
\end{enumerate}
\end{theo}

\noindent Let us now come back to the initial system and let us try to deduce
 the invariant law associated with $X_n$. If the sequence $(X_n)_n$ is defined
 by the initial terms $X_0, \hdots,  X_{k-1}$, with values in $[-L,L]$,
 and the recurrence relation  $X_{n+k}= \ph(X_n, \hdots, X_{n+k-1})$, one sets
 $Z_n= (\gamma^{j-1} X_{n+j-1})_{1\leq j \leq k}$. Then  $(Z_n)_n$ satisfies the
 recurrence relation $Z_{n+1}=T(Z_n)$, which yields the following result:

\begin{theo}
If the random variable $Z_0= (\gamma^{j-1} X_{j-1})_{1\leq j \leq k}$ has
 density  $h_*$, then, for every $n\geq 0$, $Z_n$  has
 density  $h_*$. Computing the marginal distributions, we get as a consequence
 that for every $n \in \N$, $X_n$ has a density $h_{inv}$ which has the following expressions:  for every $j \in \{0, \hdots ,k-1\}$
$$
\forall u\in [-L,L],\qquad 
h_{inv}(u)= \gamma^j  \int_{ \R^{k-1}}
h_*(z_1,\dots,\gamma^j u ,\dots,z_k)\ d\check{z}_{j+1} 
$$
where $d\check{z}_{j+1}$ means that one integrates with respect to all 
coordinates of $z$ but $z_{j+1}$. 
\end{theo}

\noindent Indeed, $\gamma^j X_n$ is the $(j+1)-$th coordinate of
 $Z_{n-j}$ if  $j= 0,\dots, k-1$. Let us consider a Borel set $A$ of $\R$. 
Then, for $j\in \{  0,\dots, k-1 \}$,
$$
\begin{array}{lll}
P(X_n \in A) & 
\displaystyle = P(Z_{n-j} \in \R^{j} \times \gamma^j A \times \R^{k-j-1})\\
& \displaystyle =
\int_{ \R^{j} \times \gamma^j A \times \R^{k-j-1}}
h_*(z_1,\dots,z_k)\ dz_1 \dots d z_k\\
& \displaystyle =
\int_{ \R^{j} \times  A \times \R^{k-j-1}}
h_*(z_1,\dots,\gamma^j u ,\dots,z_k)\ d\check{z}_{j+1}\ \gamma^{j}  du &
 {\rm with}\  z_{j+1}= \gamma^j u
 \\
& \displaystyle =
\int_{A}\left(  \int_{ \R^{k-1}}
h_*(z_1,\dots,\gamma^j u ,\dots,z_k)\ d\check{z}_{j+1}\right) \ \gamma^{j}  du,
\end{array}
$$
which gives the desired result.\\

\noindent If $F$ is defined on  $[-L,L]$ and if  $s\in\{ 1,\dots, k\}$,
 let us denote by $T_s F$ the function defined on $\Omega$ by
\begin{eqnarray}\label{TsF}
T_sF(z)=T_sF(z_1,\dots, z_k)=  F(z_s \gamma^{1-s}).
\end{eqnarray}

\noindent The following Lemma is then a direct consequence of point 6 in
 Theorem \ref{resconjug}, applied to 
  $T_sF$ and $T_rH$ for  $s,r\in\{ 1,\dots, k\}$:

\begin{lem}\label{th4}
For every Borel set $B$ of $[-L,L]$ and every interval  $I$ of $[-L,L]$,
 if $Z_0$  has the invariant distribution, one has:
$$
\begin{array}{c}
\displaystyle \left\vert P\left( X_{ n \times ppcm(L_i)+s-1}\in B,
 X_{r-1} \in I \right)  -\sum_{j=1}^{\dim(E_1)} \sum_{l=0}^{L_j-1} \mu(W_{j,l})
 <T_s \1_B,1>_{\mu_{j,l}} <1,T_r \1_I>_{\mu_{j,l}} \right| \\
\\
\leq \left( (2L)^k\gamma^{k(k-1)/2} + 4 (2L)^{k-1} \gamma^{1-r+k(k-1)/2}
\varepsilon_0^{1-\alpha} + 2^{2k}L^{k-1} 
\left( \frac{1-\gamma^k}{1-\gamma}\right)^{k-1} \varepsilon_0^{1-\alpha}\right)
 C \rho^{n }.
\end{array}
$$
More generally, let $F$, defined and measurable on $[-L,L]$, be such that 
 $T_s F$ belongs to $L^1_{\mu}(\Om)$. Let $H \in L^{\infty}_m([-L,L])$ be such
 that 
 $\displaystyle \sup_{0<\delta<\eps_0 \gamma^{1-r}}\delta^{-\alpha} \int_{[-L,L]}
 {\rm Osc}(H,]x-\delta,x+\delta[ \cap [-L,L])\ dx < +\infty$.\\
 Then $T_r H \in V_{\alpha}(\Om)$ and
$$
\left\vert E( F( X_{n \times ppcm(L_i)+s-1}) H(X_{r-1}))
 -\sum_{j=1}^{\dim(E_1)} \sum_{l=0}^{L_j-1} \mu(W_{j,l})
 \mu_{j,l}(T_s F)\mu_{j,l}(T_r H) \right| 
\leq C(F,H,L) ~ \rho^n $$
with
$$ \begin{array}{l}
\displaystyle
C(F,H,L) = C  ||T_sF||_{L^1_{\mu}}\bigg(
 (2L)^{k-1}\gamma^{k(k-1)/2}||H||_{L^1_m([-L,L])}
\\
\displaystyle
+(2L)^{k-1}\gamma^{(k(k-1)/2)-\alpha(r-1)}
 \sup_{0<\delta<\eps_0 \gamma^{1-r}}\delta^{-\alpha} \int_{[-L,L]}
 {\rm Osc}(H,]x-\delta,x+\delta[ \cap [-L,L])\ dx
\\
\displaystyle +
 2^{2k}L^{k-1} \left( \frac{1-\gamma^k}{1-\gamma}\right)^{k-1} 
\varepsilon_0^{1-\alpha} ||H||_{L^{\infty}_m([-L,L])}\bigg).
\end{array}
$$
\end{lem}

\noindent This last result, which gives the exponential decay of 
correlations, is a straightforward consequence of Lemma \ref{th4} and of
 the remark in
 point 7, Theorem \ref{resconjug}.

\begin{theo}
If, moreover, $T$ is mixing, then
$$ |Cov(F(X_{n+s-1}),H(X_{r-1}))| \leq C(F,H,L) ~ \rho^n. $$
\end{theo}

%
%

\section{Proofs}

\noindent   Theorem \ref{resconjug} is a consequence of Theorems 5.1 and 6.1
 of \cite{SAU}, which rely on \cite{ITM}, as well as \cite{HK}
in the case when $d=1$, where the use of bounded-variation functions is 
possible. The difficulty lies in verifying that $T$ satisfies Hypotheses 
 (PE1) to (PE5).\\

\noindent To prove that  (PE2) is satisfied, we shall first establish that
 $T_j$ is a $C^1$ diffeomorphism on  $W_j$ onto $T_j(W_j)$.
Hypothesis \ref{hyp3}
 about $\displaystyle \frac{\partial \ph_j}{\partial x_1}$ assures that $T_j$
is a local diffeomorphism. To check that it is injective, let us consider two different points $u$ and $v$ of  $W_j$, such that $T_j(u)=T_j(v)$. Then $u_i=v_i$ for every $2 \leq i \leq k$ and
$$ \ph_j\left(u_1,\frac{u_2}{\gamma}, \hdots, \frac{u_k}{\gamma^{k-1}}\right)= \ph_j\left(v_1,\frac{u_2}{\gamma}, \hdots, \frac{u_k}{\gamma^{k-1}}\right).$$
Using the geometrical hypothesis \ref{hyp4} and applying the fundamental theorem of calculus to $t \mapsto \ph_j(\Gamma(t))$ leads to a contradiction.
\noindent The regularity hypotheses on the  $\ph_j$ (and hence on the  $T_j$)
allow to prove that  $\det(DT_j^{-1})$ is  $\alpha$-H\"older, provided the
 domain is conveniently restricted. One can see that there exist, for each  
 $\beta_{j}>0$, an open and relatively compact set $\Val_j$ and a real constant
 $c_{j}$ such that the following holds
\begin{itemize}
\item $\overline{U_j}\subset  \Val_{j} \subset \overline{\Val_{j}} \subset W_j$ ;
\item $B_{\beta_{j}}(T_j(U_j)) \subset T_j(\Val_{j}) ;$
\item
for every $\eps<\beta_{j}$, every $z\in T_j(\Val_{j})$ and all
  $x,y\in B_{\eps}(z)\cap T_j(\Val_{j})$,
$$
\Big| \det(DT_j^{-1}(x))- \det(DT_j^{-1}(y)) \Big| \leq c_{j} \Big| \det(DT_j^{-1}(z)) \Big|\eps^{\alpha}.
$$
\end{itemize}

\noindent Setting $\beta = \min\limits_j \beta_{j} >0$ and
 $c = \max\limits_j c_{j} >0$, one obtain constants which are convenient for
 every $j \in \{1, \hdots , d\}$. Hence (PE2) is satisfied.
\\

\noindent This allows us to specify the open sets on which we shall work.
 There exists $\eps_2>0$ such that, for every $j\in\{1,\dots,d\}$,
 $B_{2\eps_2}(\overline{U_j}) \subset \Val_{j} \subset W_j$.
From now on, one sets $V_j=B_{\eps_2}(\overline{U_j})$. 
Then $T_j(V_j)$ is open and
 $T_j(\overline{U_j})$ is compact and contained in $T_j(V_j)$. One can find a
positive $\eps_{0,1} $ such that
  $B_{\eps_{0,1}}( T_j(\overline{U_j})) \subset T_j(V_j)$ for every $j$, which
proves that  Hypothesis (PE1) is satisfied.\\
\\
 Hypothesis (PE3) is clearly fulfilled since 
  $\displaystyle \Om= \bigcup_{j=1}^d  U_j\ \cup \Nal'$ is a disjoint union of
 open sets and of a negligible set.\\
\\
For (PE4), we need two steps. We first prove that the map is locally expanding
(when the preimages in  $\Val_j$ are sufficiently near, Proposition
 \ref{dilatance1}). Then we prove the hypothesis itself
 (Proposition \ref{dilatance2}), in the case when the images in $T_j(V_j)$
are sufficiently near. 

\begin{prop}\label{dilatance1}
Let  $u $ and $ v\in \Val_{j}$  be such that the segment $[u,v]$ is contained in
 $\Val_{j}$. Then
$$
|| T_j(u)-T_j(v)||^2 \geq \frac{1}{s^2} || u - v||^2.
$$
\end{prop}

\noindent{\it Proof:} One applies the fundamental theorem of
 calculus to the map defined on [0,1] by 
 $t \mapsto \ph_j(v_1 + t(u_1 - v_1), \frac{v_2+t(u_2-v_2)}{\gamma}, \hdots, \frac{v_k+t(u_k-v_k)}{\gamma^{k-1}})$, which yields a  $c\in ]0,1[$ such that  
$$
 || T_j(u)-T_j(v)||^2 = (v_1-u_1, \hdots, v_k-u_k) B \left( \begin{array}{lll}
v_1 - u_1\\
\vdots \\
v_k - u_k
\end{array}\right)
$$
where $B=(b_{i,l})_{1 \leq i,l \leq k}$ is the matrix with coefficients
\begin{eqnarray}\label{CoeffB}
\left\{
\begin{array}{llll}
b_{i,l} & = & \displaystyle \gamma^{2k-i-l} \frac{\partial \ph_j}{\partial x_i}(M_c) \frac{\partial \ph_j}{\partial x_l}(M_c) & {\rm if ~} i \neq l, \\
 & & & \\
b_{1,1} & = & \displaystyle \gamma^{2k-2} \left( \frac{\partial \ph_j}{\partial x_1}(M_c) \right)^2, & \\
 & & & \\
b_{i,i} & = & \displaystyle \frac{1}{\gamma^2} + \gamma^{2(k-i)} \left( \frac{\partial \ph_j}{\partial x_i}(M_c) \right)^2 & {\rm if ~} i > 1,
\end{array}
\right.
\end{eqnarray}
with $\displaystyle M_c = \left(v_1 + c(u_1 - v_1),
 \frac{v_2+c(u_2-v_2)}{\gamma}, \hdots, \frac{v_k+c(u_k-v_k)}{\gamma^{k-1}}\right)$.\\
The matrix $B$ is real and symmetrical. Its eigenvalues are contained in the 
Gershg\"orin disks and hence in the following domain 
$$ \bigcup_{i=1}^k \left[ b_{i,i} - \sum_{l \neq i} |b_{i,l}|, b_{i,i} + \sum_{l \neq i} |b_{i,l}| \right]. $$
We shall establish that all these intervals are contained in
 $[\sigma, +\infty[$. To that aim, it is sufficient to prove that 
$b_{i,i} - \sum_{l \neq i} |b_{i,l}| \geq \sigma$ for every $i$.\\
According to Hypothesis \ref{hyp3} one has, for every $l > 1$,
 $|\frac{\partial \ph_j}{\partial x_l}(M_c)| \leq \frac{M}{A}$, which implies
 that
$$ \left\{
\begin{array}{ll}
\displaystyle b_{1,1} - \sum_{l \neq 1} |b_{1,l}| \geq \gamma^{2k-2} A^2 - M \sum_{l \neq 1} \gamma^{2k-1-l} & \\
\displaystyle b_{i,i} - \sum_{l \neq i} |b_{i,l}| \geq \frac{1}{\gamma^2} - M \gamma^{2k-i-1} - \frac{M^2}{A^2} \sum_{l \neq i, l >1} \gamma^{2k-i-l} & {\rm for  ~} i>1.
\end{array}
\right.
$$
Since $\gamma < 1 $ and $ 2k-1-l \geq k-1$ one eventually gets:
$$ \left\{
\begin{array}{ll}
\displaystyle b_{1,1} - \sum_{l \neq 1} |b_{1,l}| \geq \gamma^{2k-2} A^2 - M (k-1) \gamma^{k-1} & \\
\displaystyle b_{i,i} - \sum_{l \neq i} |b_{i,l}| \geq \frac{1}{\gamma^2} - M \gamma^{k-1} - \frac{M^2}{A^2} (k-2) \gamma & {\rm for  ~} i>1.
\end{array}
\right.
$$
Since  $\gamma = A^{-1/k}$, one derives the inequalities:
$$ \left\{
\begin{array}{ll}
\displaystyle b_{1,1} - \sum_{l \neq 1} |b_{1,l}| \geq \frac{1}{\gamma^2} - M (k-1) \gamma^{k-1} & \\
\displaystyle b_{i,i} - \sum_{l \neq i} |b_{i,l}| \geq \frac{1}{\gamma^2} - M (k-1) \gamma^{k-1} - M^2 (k-2) \gamma^{2k+1} & {\rm for  ~} i>1.
\end{array}
\right.
$$
Therefore, the eigenvalues of $B$ are all greater than or equal to
 $\frac{1}{\gamma^2} - M (k-1) \gamma^{k-1} - M^2 (k-2) \gamma^{2k+1}$. Since $M \in ]0,M_0(\sigma,A)[$, one has:
$$ \frac{1}{\gamma^2} - M (k-1) \gamma^{k-1} -
 M^2 (k-2) \gamma^{2k+1} \geq \sigma. $$
Consequently,  the eigenvalues of $B$ are all greater than or equal to
 $\sigma = \frac{1}{s^2}$, which gives the desired result. Notice that this 
last inequality compells us to choose $\frac{1}{\gamma^2} > \sigma$.\hfill
 $\square$\\

\noindent Compacity arguments give the existence of $\eps_{0,2} > 0$ such that,
 for every $u\in \overline{V_j}$,
$$
B_{\eps_{0,2}}(T_j(u)) \subset T_j(B_{\eps_2}(u)).
$$

\begin{prop}\label{dilatance2}
Let $\eps_0= \min(\eps_{0,1},\eps_{0,2}) > 0$. Recall that $\overline{U_j} \subset V_j \subset \overline{V_j} \subset \Val_j \subset W_j$.  As a consequence,
\begin{itemize}
\item
for all $x, ~ y \in T_j(V_j)$ satisfying $\| x-y \|<\eps_0$, the following
 inequality is valid:
$$
s^2 ~ \| x-y \| > \| T_j^{-1}(x),T_j^{-1} (y) \|.
$$
\item
$B_{\eps_0}( T_j(\overline{U_j})) \subset T_j(V_j)$.
\end{itemize}
\end{prop}

\noindent{\it Proof:} The second statement comes from the fact that $\eps_0\leq \eps_{0,1}$ and from the results we obtained in relation with (PE1). \\
Let us prove the first statement, which implies Condition (PE4) of Saussol.
 Let $x, ~y \in T_j(V_j)$ satisfy $\| x-y \|<\eps_0$.
 Set $u=T_j^{-1}(x) \in V_j$. According to the preceding remark, as  $\eps_0$ is
smaller than $\eps_{0,2}$,
$$
y \in B_{\eps_0}(T_j(u)) \subset T_j(B_{\eps_2}(u)).
$$
Hence $v=T_j^{-1}(y)\in B_{\eps_2}(u)\subset \Val_j$. According to
 Proposition \ref{dilatance1}, 
$$
\| x-y \|^2 = || T_j(u)-T_j(v)||^2 > \sigma || u-v ||^2,
$$
which proves the result.\hfill $\square$\\

\noindent  To conclude, Hypothesis  (PE5)
is a consequence of Lemma 2.1 of Saussol and of Hypothesis \ref{hyp5}.\\

\noindent Since the hypotheses  (PE1) to  (PE5) are satisfied, Theorem 5.1 of
 \cite{SAU} implies the properties 1 to 5 of Theorem \ref{resconjug} about
 $V_{\alpha} $ and $ L^1_m$.
 Now, if $f \in E_i$, $f$ is equal to zero on $\Omega^c$, which implies that
 $f \in L^1_m(\Omega) $ and $ V_{\alpha}(\Omega)$.\\

\noindent To prove point 6, we shall apply Theorem 6.1 of \cite{SAU} on every
 subset $W_{j,l}$, on which a suitable power of $T$ is mixing. Adopting the notations of Point 5 of Theorem 5.2 of \cite{SAU}, there exist real constants $C>0$
and $\rho\in ]0,1[$ such that, for every $(j,l)$ satisfying 
$ 1\leq j\leq \dim(E_1)$, $0\leq l\leq L_j -1$, every function 
$f\in L^1_{\mu_{j,l}}(\Om)$ and   every function $h\in V_{\alpha}(\Om)$,
$$
\left\vert \int_{\Om}(f- \mu_{j,l}(f))\circ T^{nL_j} h ~ d\mu_{j,l}
\right\vert \leq C ||f-\mu_{j,l}(f))||_{L^1_{\mu_{j,l}}} ||h||_{\alpha,L}\rho^{n } .
$$
Let us choose, then, $h\in V_{\alpha}(\Om)$ and  $f\in L^1_{\mu}(\Om)$
 (with the result that
$f\in L^1_{\mu_{j,l}}(\Om)$ for every $j,l$). Taking the smallest common multiple
 $L'$ of the
 $L_j$ and summing the above inequalities, with $n$ replaced with
$n \frac{L'}{L_j}$, we obtain that
$$
\left| \int_{\Om} f \circ T^{ n L'}  h \ d\mu -\sum_{j=1}^{\dim(E_1)}
 \sum_{l=0}^{L_j-1} \mu(W_{j,l}) \mu_{j,l}(f)\mu_{j,l}(h) \right| 
\leq C  ||h||_{\alpha,L}||f||_{L^1_{\mu}} \rho^{n } ,
$$
remarking that $||f-\mu_{j,l}(f))||_{L^1_{\mu_{j,l}}} \leq 2 ||f||_{L^1_{\mu}}$.\\

\noindent Point 7 is a straightforward consequence of Point 6, since
 dim($E_1) = 1$ and  $L_1=1$. This concludes the proof of Theorem
\ref{resconjug}\hfill $\square$\\

\noindent Let us now turn to Lemma \ref{th4}.
 If $Z_0=(X_0,\dots,\gamma^{k-1}X_{k-1} )$ has distribution $\mu$, then this is
 the case for 
 $Z_n=(X_n,\dots,\gamma^{k-1}X_{n+k-1} )$  as well.
 If $f\in L^1_{\mu}(\Om)$ and if  $h\in V_{\alpha}(\Omega)$, one has:
$$
\left\vert E( f (Z_{nL'}) h( Z_{0}))
 -\sum_{j=1}^{\dim(E_1)} \sum_{l=0}^{L_j-1} \mu(W_{j,l}) \mu_{j,l}(f)\mu_{j,l}(h)
 \right| \leq C ||f||_{L^1_{\mu}} ||h||_{\alpha,L}\rho^{n } .
$$
\noindent
Let  $r,s$ be in $\{1,\dots,k\}$ and let  $F,H$ be measurable functions
 defined on   $[-L,L]$.
The application  $T_r H$ (defined in (\ref{TsF})) belongs to
 $V_{\alpha}(\Omega)$ if and only if 
$H$ is in $L^{\infty}([-L,L],m)$ and satisfies
$$
 \sup_{0<\delta<\eps_0 \gamma^{1-r}}\delta^{-\alpha} \int_{[-L,L]}
 {\rm Osc}(H,]x-\delta,x+\delta[ \cap [-L,L])\ dx < +\infty.
$$
One then has
$$
\begin{array}{lll}
|| T_r H||_{\alpha,L}     \displaystyle =
 (2L)^{k-1}\gamma^{k(k-1)/2}||H||_{L^1_m([-L,L])}\\
+\displaystyle
(2L)^{k-1}\gamma^{(k(k-1)/2)-\alpha(r-1)}
 \sup_{0<\delta<\eps_0 \gamma^{1-r}}\delta^{-\alpha} \int_{[-L,L]}
 {\rm Osc}(H,]x-\delta,x+\delta[ \cap [-L,L])\ dx
\\+ \displaystyle
2^{2k}L^{k-1} \left( \frac{1-\gamma^k}{1-\gamma}\right)^{k-1} \varepsilon_0^{1-\alpha} ||H||_{L^{\infty}_m([-L,L])}.
\end{array}
$$

\noindent Consequently, if  $H$ satisfies these conditions and if $F$ is
 such that $T_sF$ belongs to $L^1_{\mu}(\Om)$, for example if $F$ is measurable
 and bounded on $[-L,L]$, one gets the second statement of
Lemma \ref{th4}

$$
\begin{array}{lllll}
\displaystyle
\left\vert E( F( X_{n \times L'+s-1}) H(X_{r-1}))
 -\sum_{j=1}^{\dim(E_1)} \sum_{l=0}^{L_j-1} \mu(W_{j,l})
 \mu_{j,l}(T_s F)\mu_{j,l}(T_r H) \right| \leq
 C  ||T_sF||_{L^1_{\mu}} \\
\displaystyle \bigg(
 (2L)^{k-1}\gamma^{k(k-1)/2}||H||_{L^1_m([-L,L])}
\\
\displaystyle
+(2L)^{k-1}\gamma^{(k(k-1)/2)-\alpha(r-1)}
 \sup_{0<\delta<\eps_0 \gamma^{1-r}}\delta^{-\alpha} \int_{[-L,L]}
 {\rm Osc}(H,]x-\delta,x+\delta[ \cap [-L,L])\ dx
\\
\displaystyle
+
 2^{2k}L^{k-1} \left( \frac{1-\gamma^k}{1-\gamma}\right)^{k-1} \varepsilon_0^{1-\alpha} ||H||_{L^{\infty}_m([-L,L])}\bigg)\rho^n.
\end{array}
$$

\noindent In particular, if $H$ is the indicator function of an interval
and $F$, that of a Borel set,  we obtain the first assertion of Lemma 
 \ref{th4}.

%
%

\section{A nonlinear example}

\noindent Set $\Omega= [-L,L]^k$.  We can state the result:
\begin{lem}\label{premier}
Let $A>1$ and $M>0$ be real numbers.
\noindent Let $\alpha_0, a_1,\dots a_k, b_1$ be  real positive numbers, with 
 $a_1>0$.
Set
$$
\psi(x)= \alpha_0 + \left( \sum_{i=1}^k a_i x_i^2 \right) + b_1x_1.
$$
If the following conditions are fullfilled, 
\begin{eqnarray}
4\alpha_0a_1= b_1^2 ;\label{rel0}\\
-2a_1L+b_1 
\geq 2 A \sqrt{\frac{1}{4 a_1}(-2a_1L+ b_1)^2 + \sum_{i=2}^k a_i L^2}>0 ;
\label{rel1}\\
\forall i \in [2,k]\cap \N,\ \sqrt{a_1}\sqrt{a_i}\leq 2M,\label{rel2}
\end{eqnarray}
 the application $\psi$ is positive on  $\Omega$ and
 $\ph_{0}=\sqrt{\psi}$ is well defined and  $C^{\infty}$ on an open neighbourhood 
of $\Omega$. Moreover  $\ph_{0}=\sqrt{\psi}$  satisfies the inequalities 
\begin{equation}
\forall x\in \Omega,\ \left| \frac{\partial\ph_{0} }{\partial x_1}(x)\right|\geq A 
,\quad \forall i\in [2,k]\cap\N,\ \left|
\frac{\partial\ph_{0} }{\partial x_1}(x)
 \frac{\partial\ph_{0} }{\partial x_i}(x)\right| \leq M. 
\label{derivees-partielles}
\end{equation}
\end{lem}

\noindent{\it Proof:}
Thanks to (\ref{rel0}), one gets
$\psi(x)= \frac{1}{4a_1}(2a_1x_1+b_1)^2 +\sum_{i=2}^k a_i x_i^2$. Hence
$\psi(x)\geq \frac{1}{4a_1}(2a_1x_1+b_1)^2$ and, since
$2a_1x_1+b_1\geq -2a_1L+b_1>0$ sur $[-L,L]$ according to (\ref{rel1}), 
$\psi$ is positive on $\Omega$.
Therefore $\ph_{0}=\sqrt{\psi}$ is well defined and 
 $C^{\infty}$ on an open neighbourhood 
of $\Omega$.
One checks that
$$
\frac{\partial \ph_0}{\partial x_1}(x)= \frac{2a_1x_1+b_1}{2\sqrt{\psi(x)}}
\geq \frac{2a_1x_1+b_1}{2\sqrt{\psi(x_1,L,\dots,L)}}.
$$
Denoting by $g= g(x_1)$ the function appearing in the right side above,
 one sees that 
$g'$ has the sign of $2a_1\sum_{i=2}^k a_i L^2$, which means that $g$ is
an increasing function. To obtain the desired condition about
$\frac{\partial \ph_0}{\partial x_1}$, it suffices that $g(-L)\geq A$,
which is a consequence of (\ref{rel1}).\\
One has
$$
\left|
\frac{\partial\ph_{0} }{\partial x_1}(x)
 \frac{\partial\ph_{0} }{\partial x_i}(x)\right|
 =\frac{a_i|x_i|(2a_1x_1+b_1)}{2\psi(x)}
\leq
\frac{a_i|x_i|(2a_1x_1+b_1)}{2\psi(x_1,0,\dots,0,x_i,0,\dots,0)}.
$$
This can be written as
 $$
\left|
\frac{\partial\ph_{0} }{\partial x_1}(x)
 \frac{\partial\ph_{0} }{\partial x_i}(x)\right|
\leq
\sqrt{a_ia_1}
\frac{ (\sqrt{a_i} |x_i|) (\sqrt{a_1}x_1+\frac{b_1}{2\sqrt{a_1}})}
{ (\sqrt{a_1}x_1+\frac{b_1}{2\sqrt{a_1}})^2+  ((\sqrt{a_1}x_1)^2},
$$
and it is easy to see that it is smaller than 
 $\frac{\sqrt{a_ia_1}}{2}$, hence smaller than $M$ according to (\ref{rel2}).

\begin{lem}\label{deuxieme}
Let $A>1$ and $M>0$ be real numbers. Let $a_1 $ and $ b_1$ be such that 
\begin{eqnarray}
a_1 \geq 2A^2,\label{a_1}\\
b_1\geq 4LM \sqrt{k-1} + 2a_1L.\label{b_1}
\end{eqnarray}
If we set $4\alpha_0a_1= b_1^2 $ (\ref{rel0}) and impose that,
 for all $i$ between $2$ and $k$, $\sqrt{a_1}\sqrt{a_i}\leq 2M$ (\ref{rel2}),
the system of Lemma \ref{premier} is satisfied.
\end{lem}

\noindent{\it Proof}\\
According to (\ref{b_1}), $b_1-2a_1L \geq 4LM \sqrt{k-1}>0$ and then
$\displaystyle \sqrt{\frac{1}{4 a_1}(-2a_1L+ b_1)^2 +\sum_{i=2}^k a_i L^2}>0$.
Moreover, using  (\ref{rel2}), one gets: 
$$
2A \sqrt{ \frac{1}{4 a_1}(-2a_1L+ b_1)^2 + \sum_{i=2}^k a_i L^2  }
\leq 
2A \sqrt{ \frac{1}{4 a_1}(-2a_1L+ b_1)^2 + (k-1) L^2 \frac{4M^2}{a_1} }.
$$
This last expression is smaller than $-2a_1L+b_1$ if and only if
$$
(-2a_1L+ b_1)^2 \geq  \frac{4A^2}{4 a_1}(-2a_1L+ b_1)^2 +
 (k-1) 4A^2 L^2 \frac{4M^2}{a_1},
$$
or equivalently
$$
a_1(-2a_1L+ b_1)^2 \geq  A^2(-2a_1L+ b_1)^2 +
 (k-1) 16A^2M^2 L^2 .
$$
But according to (\ref{a_1}) and (\ref{b_1}),
$$
\begin{array}{lll}
\displaystyle
a_1(-2a_1L+ b_1)^2 - A^2(-2a_1L+ b_1)^2 -
 (k-1) 16A^2M^2 L^2 \\
\qquad \geq \displaystyle
 2A^2(-2a_1L+ b_1)^2 - A^2(-2a_1L+ b_1)^2 -
 (k-1) 16A^2M^2 L^2 \\
\qquad \geq \displaystyle
 A^2\left( (-2a_1L+ b_1)^2  -
 (k-1) 16M^2 L^2\right) \geq 0.\\
\end{array}
$$

\vskip 1 cm

\noindent We assume that the conditions of Lemma \ref{deuxieme} are satisfied. 

\noindent  For $\ell\in [-L,L[$ and $p\in \Z$, set
$$
\ph_{\ell,p}(x)= \ell + \sqrt{\psi(x)} + 2pL.
$$
One defines the application $\ph$ almost everywhere on $\Omega$, by
$$
\ph(x) = \ph_{\ell,-p}(x) \quad {\rm si}\quad \ph_{\ell,0}(x) \in [2pL-L,2pL +L[.
$$ 
Since $\ph_{\ell,-p}(x) =\ph_{0}$ up to an additive constant, the conditions (\ref{derivees-partielles}) concerning the partial derivatives  are fullfilled. \\

\noindent Let us specify the open sets. For $p\in\Z$, set
$$
\ooo_p=  \{ x\in ]-L,L[^k : \ph_{\ell,0}(x) \in ]2pL-L,2pL +L[\}.
$$
One sees that, for $p\leq -1$, $\ooo_p$ is empty and that, otherwise, 
$$
\begin{array}{lllll}
\ooo_0= \{ x\in ]-L,L[^k : \psi(x) < (L-\ell)^2\},\\
\ooo_p=  \{ x\in ]-L,L[^k : ((2p-1)L- \ell)^2 <\psi(x) < ((2p+1)L-\ell)^2\},
& p\geq 1.\\
\end{array}
$$
The sets $\ooo_p$ are open and may be empty. We would specify which ones are
 empty if we really wanted to give an 
explicit expression of the Frobenius-Perron operator.\\
Put $\sss_p= \{x\in \R^k : \psi(x) =  ((2p-1)L- \ell)^2\}.$

\noindent If  $\sss_p\cap\Omega$ is not empty, 
$\displaystyle \frac{\partial \psi}{\partial x_1}(x) >0$ is valid 
for every point of $\sss_p\cap\Omega$ according to (\ref{rel1}), so
 $x_1$ can be considered, locally, as a   $C^{\infty}$ 
function of the other $x_i$ and $\sss_p\cap\Omega$  is a finite union of
  $C^{\infty}$ submanifolds. 
The edges of $\Omega$ are parts of hyperplanes and so are   $C^{\infty}$ too.
A submanifold  $\sss_p$ crosses at most  $k$ hyperplanes, which implies that
the maximal crossing number, $Y$, is smaller than $ k+1$.\\

\noindent The geometrical condition, under its simple form, is satisfied.
Indeed, let $U=(u_1,u_2,\dots, u_k) $ and $ V=(v_1,u_2,\dots, u_k)$ be two 
points of the same set $\ooo_p$.  On  $[-L,L]^k$,
 $\frac{\partial \psi}{\partial x_1}(x) >0$, 
according to (\ref{rel1}).
Hence, for $t\in [0,1]$, if one assumes that $-L<u_1<v_1<L$, 
$$
\psi(U) \leq \psi(tU+(1-t)V) \leq \psi(V),
$$
since the only coordinate that changes is the first one. Therefore
 $\psi(tU+(1-t)V)$ is in the same interval 
($] ((2p-1)L- \ell)^2, ((2p+1)L-\ell)^2[$ if  $p\geq 1$) as $\psi(U) $ 
and $ \psi(V) $. Consequently, $ tU + (1-t)V$ is in  $\ooo_p$.

%
%

\end{document}